\documentclass[twoside]{article}

\usepackage{PRIMEarxiv}

\usepackage[utf8]{inputenc} 
\usepackage[T1]{fontenc}    
\usepackage{hyperref}       
\usepackage{url}            
\usepackage{booktabs}       
\usepackage{amsfonts}       
\usepackage{nicefrac}       
\usepackage{microtype}      
\usepackage{fancyhdr}       
\usepackage{graphicx}       
\usepackage{rotating}
\usepackage{longtable, array}
\usepackage{tabularx}
\usepackage{pdflscape}
\usepackage{amssymb,amsthm,amsmath}
\usepackage{adjustbox}
\usepackage{placeins}
\title{A Comparison of Image and Scalar-Based Approaches in Preconditioner Selection
}

\author{
  Michael Souza \\
  Universidade Federal do Ceará (UFC) \\
  Fortaleza, Ceará, Brazil \\
  \texttt{michael@ufc.br} \\
   \And
  Luiz Mariano Carvalho \\
  Universidade do Estado do Rio de Janeiro (UERJ) \\
  Rio de Janeiro, Brazil \\
  \texttt{luizmc@ime.uerj.br} \\
  \AND
   Douglas Augusto \\
  Fundação Oswaldo Cruz (Fiocruz) \\
  Rio de Janeiro, Brazil \\
  \texttt{daa@fiocruz.br} \\
   \And
  Jairo Panetta\\
  Instituto Tecnológico da Aeronáutica (ITA) \\
  São José dos Campos, São Paulo, Brazil \\
  \texttt{jairo.panetta@gmail.com} \\
  \AND  
  Paulo Goldfeld \\
  Universidade Federal do Rio de Janeiro (UFRJ) \\
  Rio de Janeiro, Brazil  \\
  \texttt{goldfeld@matematica.ufrj.br} \\
   \And
  José Roberto Pereira Rodrigues \\
  Petrobras R\&D Center (CENPES) \\
  Rio de Janeiro, Brazil \\
  \texttt{jrprodrigues@petrobras.com.br}
}

\begin{document}
\maketitle

\begin{abstract}
Within high-performance computing (HPC), solving large sparse linear systems efficiently remains paramount, with iterative methods being the predominant choice. However, the performance of these methods is tightly coupled to the aptness of the chosen preconditioner. The multifaceted nature of sparse matrices makes the universal prescription of preconditioners elusive. Notably, the key attribute of sparsity is not precisely captured by scalar metrics such as bandwidth or matrix dimensions. Advancing prior methodologies, this research introduces matrix sparsity depiction via RGB images. Utilizing a convolutional neural network (CNN), the task of preconditioner selection turns into a multi-class classification problem. Extensive tests on 126 SuiteSparse matrices emphasize the enhanced prowess of the CNN model, noting a 32\% boost in accuracy and a 25\% reduction in computational slowdown.
\end{abstract}


\section{Introduction}\label{sec:intro}
Numerical simulations in reservoir engineering require the solution of large linear systems with sparse matrices, taking up more than 50\% of the computing time \cite{gasparini21HybridParallelIterativeSparseLinearSolverFramework,gaganis2012machine}. Krylov methods are preferred linear solvers in this context. However, their success depends on choosing a suitable preconditioner, which is still a challenging task without an established standard, often based on trial and error or user experience \cite{saad2003iterative,scott2023introduction}.

A preconditioner is a matrix or operator that modifies the original system to speed up the convergence of iterative linear solvers. Multiple factors influence the preconditioner effectiveness, including features of the sparse matrix, the computational architecture, and the data structures employed \cite{benzi2002preconditioning,bell2008efficient}. 
In this work, we will address only the first of these factors.

Designing an appropriate representation for sparse matrices in the context of preconditioner suggestion is challenging, mainly due to scalability requirements. For large sparse systems, the computational complexity should ideally be linear concerning the number of non-zero elements. This restriction arises from the need to maintain efficiency as the system size increases, a common scenario in reservoir simulations or computational fluid dynamics \cite{stuben2001AReview}. The difficulty in obtaining compact representations lies in capturing the matrix properties that most influence the performance of preconditioners within this linear complexity restriction.

Among the matrix attributes that influence the choice of preconditioners, we can highlight the order of the matrix, its eigenvalues, conditioning number, sparsity pattern, density, and diagonal dominance, among others. These attributes are numeric; the exception is the sparsity pattern that encapsulates topological attributes about the connection between the non-zero elements of the matrix. Although scalar attributes such as bandwidth and number of non-zero elements provide information about complex sparsity patterns, their sensitivity is low. This characteristic makes it difficult to obtain an adequate numerical representation for sparsity. 

Although not the only determining factor, the sparsity pattern impacts the level of parallelism achievable both in the application and in the construction of preconditioners such as ILU($k$) \cite{Meijerink1977AnIS, saad2003iterative}. In Algebraic Multigrid (AMG) methods, coarse operators encode sparsity and numerical values to create a multilevel approximation of the \cite{stuben2001AReview} system. Depending on the sparsity pattern, even direct methods may be applicable to solve sparse problems \cite{Davis2016ASO}.

At this stage, we explore machine learning (ML) techniques to automatically obtain compact representations of sparse matrices to select preconditioners. Extending the approach of Yamada et al. \cite{yamada2018preconditioner}, we use RGB images to spatially encode sparsity patterns, facilitating data-driven analysis to increase the efficiency of linear solvers. Unlike Yamada's work, our methodology incorporates multi-label classification to identify a range of appropriate preconditioners for a given sparse matrix. To process the image-based representations, we employ a Convolutional Neural Network (CNN) \cite{li2022survey}.

Our contributions are as follows: 
 \begin{description}  
\item [Multi-Label Model:] In conventional classification tasks, a sparse matrix is often mapped to just one preconditioner. However, in our dataset, approximately 20\% of matrices exhibit multiple optimal preconditioners. Based on research by Yamada et al. \cite{yamada2018preconditioner}, we deviate from this one-to-one mapping approach. By adopting a multi-label structure, we allow assigning multiple preconditioners to a single matrix, better handling cases where multiple preconditioners perform similarly.  
\item [Scalar models versus image-based models:] We compared scalar and image-based models. Our investigation also introduced a mixed model that combines both attributes by transforming the image-based attributes into a vector format (\emph{flattening}) and subsequently integrating them into the scalar attribute table.  
\item[Promising results:] In our initial research, image-based models outperform scalar-based ones in preconditioner selection. They show a 32\% higher probability of optimal preconditioner suggestion and a 26\% higher chance of low computational impact when one does not choose the optimal preconditioner. These results highlight the superior efficiency and effectiveness of image-based models.  
\end{description}

The remainder of this article is structured as follows: Section \ref{sec:related} reviews relevant literature, Section \ref{sec:method} outlines our methodologies, Section \ref{sec:results} discusses empirical results, and Section \ref{sec:conclusions} concludes the paper, highlighting its significance and suggesting future research directions.

\section{Related Works}\label{sec:related}

Machine Learning's potential to discern intricate patterns and facilitate data-driven decision-making has been recognized as an effective solution to various challenges in High-Performance Computing applied to the solution of linear systems  \cite{falch2017machine,tuncer2017diagnosing,memeti2019using}.

\subsection*{Sparse Matrix Storage Optimization} 
A significant application area of ML is automatic data structure selection for sparse matrix storage. Sedaghati et al. \cite{sedaghati2015automatic} used decision tree algorithms to automate storage format selection based on array properties. Nisa et al. \cite{nisa2018effective} applied machine learning techniques to predict the most suitable storage formats for GPUs. The importance of synchronizing storage data structure with computational efficiency is also highlighted in the research of Barreda et al. \cite{barreda2020performance} and Cui et al. \cite{cui2016code} that explore performance improvements on different computing platforms.

\subsection*{Auto-tuning Linear Solvers}
Another notable ML application involves (\emph{autotuning}) linear solvers. For example, Peairs and Chen \cite{peairs2011using} leveraged a reinforcement learning strategy to determine the optimal restart parameters for the GMRES iterative method. Bhowmick et al. \cite{bhowmick2006application} applied ML to select the best solvers for sparse linear systems at runtime, adapting to the data and the available computational architecture. Dufrechou et al. \cite{dufrechou2019automatic} employed ML techniques to predict the optimal solver for a specific linear system, focusing on situations when a limited number of triangular systems are solved for the same matrix. In another approach, Funk et al. \cite{funk2022prediction} presented a deep learning technique to identify the optimal iterative solver for a given sparse linear system, achieving a top-1 classification accuracy of 60\%.

\subsection*{Neural Network Approaches} 
A growing trend in ML is using neural networks to accelerate linear algebra applications. Cui et al. \cite{cui2016code} employed an ML system to predict the best implementation of matrix-vector multiplication (SpMV) for a given matrix. G{\"o}tz and Anzt \cite{gotz2018machine} introduced a convolutional neural network (CNN) to detect block structures in matrix sparsity patterns. In a different approach, Ackmann et al. \cite{ackmann2020machine} proposed using a feed-forward neural network as a preconditioner. Taghibakhshi et al. \cite{taghibakhshi2021optimization} introduced an innovative method using a reinforcement learning agent based on graph neural networks (GNNs) to construct coarse spaces in a multigrid approach. 

Although ML has substantially optimized several aspects of computational linear algebra, its potential in preconditioner selection warrants further exploration. This work is a contribution in this direction.

\section{Methodology}\label{sec:method}
We propose encoding matrix sparsity and some easily computable attributes as RGB images. This geometric representation is a compact and adequate description of the underlying matrix structure. To substantiate this claim, we conducted experiments comparing ML models trained on scalar matrix attributes with those trained on images. In both scenarios, the main objective is the optimal selection of preconditioners to minimize the convergence time of linear solvers.

\subsection{Matrix Data Set and Scalar-based Features}\label{sec:data_set}
For this study, we employed 126 symmetric, non-diagonal, and positive definite matrices obtained from the SuiteSparse Matrix Collection using the \texttt{ssgetpy} library \cite{kolodziej2019suitesparse,ssgetpy}. In the Appendix Matrix Data Set, we describe these matrices.

In this study, we employed a set of scalar features to characterize matrices, as delineated in Table \ref{tab:scalar_features}. Among these features, the condition number was calculated using MATLAB's \texttt{condest} function, while the smallest and largest eigenvalues were obtained through the \texttt{eigs} function with \texttt{smallestabs} and \texttt{largestabs} options, respectively \cite{matlab}. While the condition number and eigenvalues are informative, they are computationally expensive. In contrast, our image-based approach provides a similarly informative yet more efficient alternative.

\begin{table}
    \centering
    \caption{Scalar Features.}
    \label{tab:scalar_features}
\begin{tabular}{@{}p{2cm}p{9cm}@{}}
    \toprule
    \textbf{Features} &  \multicolumn{1}{c}{\textbf{Definition}} \\
    \midrule
    Density & Number of non-zero elements divided by the total number of elements. \\
    N& Number of rows. \\
    NNZ & Number of nonzero elements. \\
    RowNNZ & Average number of nonzero elements per row. \\
    Condest & Estimate of the 1-norm condition number of the matrix. \\
    MinEigs & Estimate of the smallest eigenvalue of the matrix. \\
    MaxEigs & Estimate of the largest eigenvalue of the matrix. \\
    DDom & Percentage of the rows diagonally dominated. If $D$ is the diagonal matrix extracted from $A$, $B = A - D$, $S = \{ i: |D_{ii}| > \sum_j{|B_{ij}|}\}$, then $\text{DDom} = \frac{|S|}{N}.$\\
    DDeg & Minimal ratio between the absolute value of the diagonal element and the sum of the non-diagonal entries. $\text{DDeg} = \min_{i}\frac{|D_{ii}|}{\sum_j{|B_{ij}|}}$.\\
     \bottomrule
\end{tabular}
\end{table}

Table \ref{tab:scalar_features_summary} presents the statistical summary of the scalar features for the matrix data set. Notably, despite all matrices being symmetric and positive definite (SPD), the smallest estimated eigenvalue is indicated as $-3.60$. This contradiction stems from the limited precision of the \texttt{eigs} function, which may not accurately compute the smallest eigenvalue for specific matrices. Although the condition number of an SPD matrix is related to its eigenvalues, the attributes we employ are estimates, making them non-redundant.

\begin{table}
    \caption{Matrices Features.}
    \label{tab:scalar_features_summary}
    \begin{adjustbox}{width=\textwidth}
    \begin{tabular}{lrrrrrrrrr}
    \toprule
     & N & NNZ & Density & RowNNZ& MaxEigs & MinEigs & Condest & DDom & DDeg \\
    \midrule
    \textbf{mean} & 4.1E+04& 6.1E+05& 1.1E-02& 2.5E+01& 1.9E+13& 9.6E+07& 1.2E+17& 1.4E-01& 2.0E-01\\
    \textbf{std} & 1.2E+05& 1.2E+06& 1.8E-02& 3.8E+01& 1.3E+14& 7.3E+08& 9.1E+17& 3.1E-01& 3.6E-01\\
    \textbf{min} & 1.0E+02& 5.9E+02& 5.0E-06& 2.9E+00& 9.8E-13& -3.6E+00& 3.6E+00& 0.00& 7.9E-15\\
    \textbf{25\%} & 1.1E+03& 1.9E+04& 3.8E-04& 6.8E+00& 9.7E+00& 1.5E-04& 4.5E+03& 0.00& 8.5E-04\\
    \textbf{50\%} & 8.6E+03& 1.3E+05& 3.8E-03& 1.6E+01& 2.9E+04& 2.0E-01& 3.9E+06& 0.00& 4.8E-03\\
    \textbf{75\%} & 2.9E+04& 5.6E+05& 1.3E-02& 2.6E+01& 5.7E+08& 2.8E+00& 1.5E+09& 2.6E-03& 2.2E-01\\
    \textbf{max} & 1.0E+06& 5.5E+06& 9.1E-02& 3.6E+02& 1.1E+15& 5.8E+09& 9.4E+18& 1.0E+00& 1.8E+00\\
    \bottomrule
    \end{tabular}
    \end{adjustbox}
\end{table}

\subsection{Image-based Features}

To visually represent matrix features, we constructed an image data set based on the method proposed by Yamada et al. \cite{yamada2018preconditioner}. Each matrix \( A \) is partitioned into blocks \( A_{ij} \) of dimension \( b \approx N / m \), where \( N \) is the dimension of the matrix and $m$ is the image resolution. These blocks are represented as pixels \( p_{ij} \) in an image of \( m \times m \) pixels. We generated four sets of images with \( m \in \{32, 64, 128, 256\}\).

Each pixel \( p_{ij} \) consists of three components aligned with the RGB channels: red, green, and blue. The red channel captures the magnitude of non-zero elements within the matrix block, the blue embodies the matrix's dimensions, and the green conveys the block's density. As a result, distinctive characteristics of each matrix are visually depicted, with the attributes of each block influencing its corresponding pixel's color. 

We will denote the red, green, and blue channels of pixel $p_{ij}$ as $R_{ij}$, $G_{ij}$, and $B_{ij}$, respectively. The range of the pixel components is limited to 0 up to 255. This restriction and the need for a general representation of heterogeneous matrices require a careful encoding scheme. We now describe the transformation of the matrix into an image:

\begin{description}
	\item[Blue Channel:] The blue channel has the same value for all pixels, given by
	$$B_{ij} = \left\lfloor\frac{N - N_{\min}}{N_{\max}  - N_{\min}}\times 255\right\rfloor,~i,j=1,\ldots,m,$$
	where $N$ is the matrix order, $N_{\min}$ and $N_{\max}$ are the minimum and maximum orders of all matrices in the data set, respectively.
	
	\item[Green Channel:] The green channel's value, \( G_{ij} \), is the density of non-zero elements within a given matrix block. It is defined as:
\[ G_{ij} = \left\lfloor\frac{NNZ_{ij}}{b^{2}}\times 255\right\rfloor \]
where \( \operatorname{NNZ}_{ij} \) represents the count of non-zero elements in the block \( A_{ij} \) and \( b \) denotes the block's order.

	\item[Red Channel:] To represent the average magnitude of non-zero elements within blocks of a sparse matrix, we first adjust or ``bias'' its non-zero elements.  This is done to ensure that all values are greater than zero. Specifically, for any given non-zero element \(a\) in the matrix \(A\), its biased value \(v(a)\) is given by $v(a) = a - \min(A) + 1,$ where $\min(A)$ is the minimal element of $A$.

The next step is to compute the average of these biased values for each matrix block. Let's call this average \( \gamma_{ij} \) for the block  \(A_{ij}\). If the difference between the maximum and minimum values of the entire matrix is 255 or less, then we compute the average of the biased values within the block. If the difference is more than 255, we take the base-2 logarithm of the biased values before averaging. This method efficiently handles blocks of widely varying values by shifting focus from the magnitude to the order of magnitude.

The formula is:
   \[\gamma_{ij} = \frac{\sum_{a \in A_{ij}} v(a)}{NNZ_{ij}},\]
   when the difference between maximum and minimum values is 255 or less, and 
   \[\gamma_{ij} = \frac{\sum_{a \in A_{ij}} \log_2 v(a)}{NNZ_{ij}},\]
   otherwise.

Finally, with the block averages at hand, we can define the value $R_{ij}$ as the normalization of the block average concerning the overall range of block averages, i.e. 
   \[R_{ij} = \left\lfloor \frac{\gamma_{ij} - \min(\gamma)} {\max(\gamma) - \min(\gamma)} \times 255 \right\rfloor\]
\end{description}

\begin{figure}[ht]
    \centering
    \includegraphics[width=0.8\textwidth]{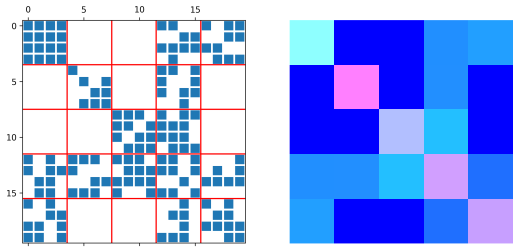}
    \caption{Image-based representation of a $20 \times 20$ matrix. The matrix is divided into $5 \times 5$ blocks, and the feature values are coded into the RGB color channels.}
    \label{fig:matrix_to_image}
\end{figure}

Figure \ref{fig:matrix_to_image} illustrates the conversion of a $20 \times 20$ matrix into a $5\times 5$ pixel image, with \(m=5\). Here, each pixel in the image corresponds to a $4 \times 4$ block of the matrix.
The full image data set can be accessed and downloaded from the GitHub repository \cite{ImagePrecGitHub}.

\subsection{Extended Scalar-based Features}\label{subsec:extended}

Representing sparsity as images naturally captures topological information, encoding it into geometric relationships, which are more challenging to achieve with scalar representations. To highlight the importance of these geometric nuances in image representations, we generated a mixed database combining scalar attributes and images. In this variation, the RGB images of $m \times m$ pixels were ``stretched'' into vectors with $3m^2$ entries (\emph{flattening}). We concatenate the flattened vectors to the scalar attributes, thus combining the two representations.

This process generated four sets of extended scalar attributes, corresponding to each $m \in \{32, 64, 128, 256\}$. To manage the increase in dimensionality and retain the most informative aspects of these extended attributes, we use Principal Component Analysis (PCA), capturing 99\% of the variation in the data \cite{guyon2006introduction}.

\subsection{Preconditioner Data Set}\label{sec:preconditioners}
We considered the following ten preconditioners: Incomplete LU Factorization with level 0 (ILU-0) and 1 (ILU-1), Block Jacobi (BJACOBI), Successive Over-Relaxation (SOR), Point Block Jacobi (PBJACOBI), Multigrid (MG), Jacobi (JACOBI), Incomplete Cholesky Factorization (ICC), Geometric Algebraic Multigrid (GAMG), Eisenstat (EISENSTAT) \cite{saad2003iterative}.

The preconditioners came from those readily available in PETSc, a prevalent software package for solving partial differential equations \cite{petsc-web-page}. In PETSc, these preconditioners were accessed by employing the \texttt{-pc} option, with the following specified values: \texttt{ilu -pc\_factor\_levels 0}, \texttt{ilu -pc\_factor\_levels 1}, \texttt{bjacobi}, \texttt{sor}, \texttt{pbjacobi}, \texttt{mg},  \texttt{jacobi}, \texttt{icc}, \texttt{gamg}, and \texttt{eisenstat}.

\subsection{Optimal Precontidioners Data Set}

To create the optimal preconditioner dataset, we generated ten linear systems (i.e, ten different right-hand side vectors) for each matrix: random solution vectors $x^\star$ were created sampling each entry independently from the standard normal distribution, and the correspondent right-hand side vector computed as $b=Ax^\star$.
We then solved each system using the Preconditioned Conjugate Gradient (CG) method, an iterative solver for Symmetric Positive Definite (SPD) \cite{saad2003iterative} linear systems.

For each linear system and preconditioner chosen, we ran the CG method five times and calculated the median resolution time. Our performance measure was obtained by summing the median times of all ten systems.

Some configurations did not converge to acceptable solutions in a timely manner. To make performance measures compatible, we use relative residual and error.

Let \( \bar{x} \) be the estimated solution vector for the system \( Ax = b \). In this case, relative residual is defined as $$ r(\bar{x}) = \frac{\|A\bar{x} - b\|_2}{\|b\|_2}$$ and the relative error is defined as $$ e (\bar{x}) = \frac{\|\bar{x} - x^\star\|_2}{\|x^\star\|_2}, x^\star \neq 0.$$

We define a pair (matrix, preconditioner) as \textit{feasible} if, in all ten linear systems, we have $r(\bar{x}) < 0.01$, $e(\bar{x}) < 0 .1$ and the total time to obtain solutions is less than one minute. Each pair that did not meet these conditions was labeled \textit{infeasible}. With this time limit (one minute), all matrices were feasible when paired with at least one of the preconditioners.

For each matrix $A$ and preconditioner $P$, we will say that $P$ is an \textit{optimal preconditioner} if $(A,P)$ is feasible and its execution time is at most 10\% slower than the fastest viable preconditioner for the same matrix.

This definition resulted in 126 sets of optimal preconditioners, one set per matrix. The optimal preconditioner dataset is in \cite{ImagePrecGitHub}.

\begin{figure}
    \centering
    \includegraphics[width=0.60\linewidth, trim={0 0 0 0}, clip]{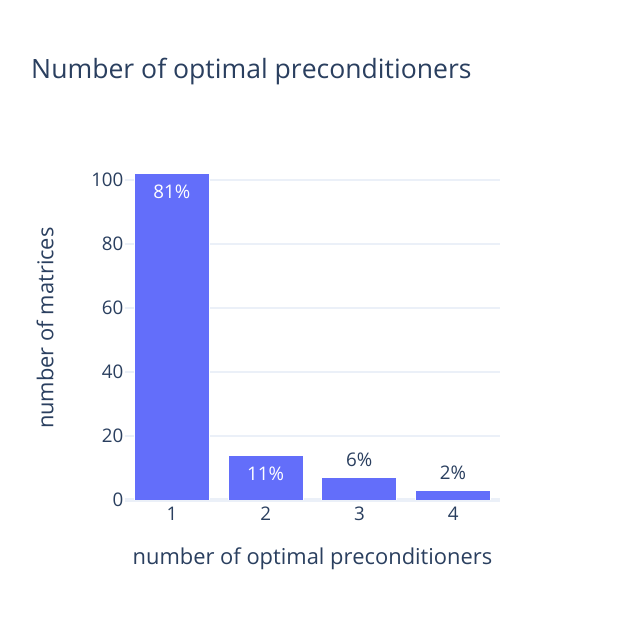}
    \caption{Distribution of matrices by number of optimal preconditioners.}
    \label{fig:pc_len}
\end{figure}

Figure~\ref{fig:pc_len} shows the distribution of matrices based on the number of optimal preconditioners. The majority of matrices, 81\%, have a single optimal preconditioner. Fewer than 11\% of the matrices have two optimal preconditioners, while 6\% and 2\% have three and four optimal preconditioners, respectively. Approximately 20\% of the matrices in the data set exhibit multiple optimal preconditioners and the average number of optimal preconditioners per matrix is 1.29. This observation accentuates the need for a multi-label formulation.

\begin{figure}
    \centering
    \includegraphics[width=0.9\linewidth]{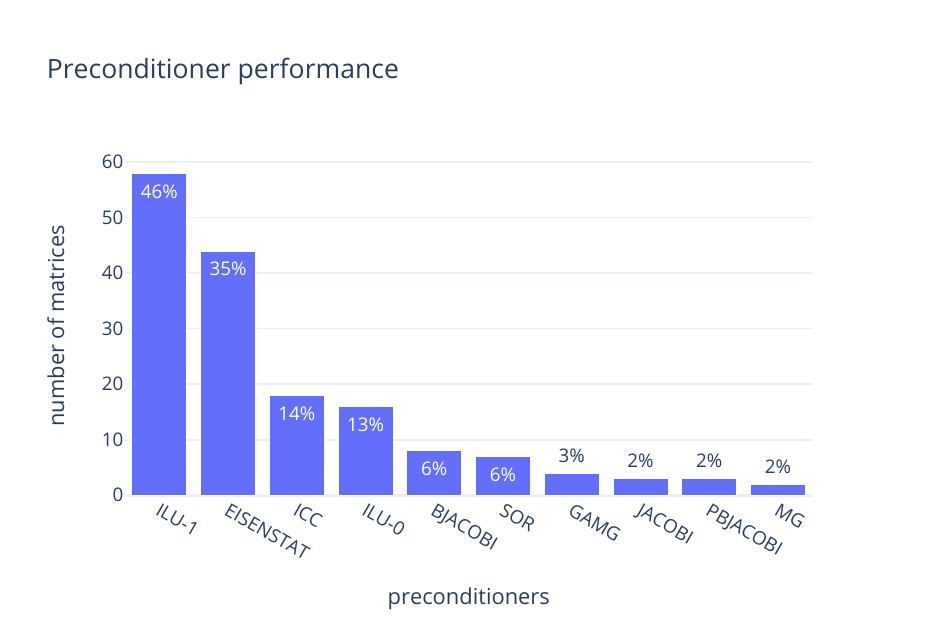}
    \caption{Distribution of matrices based on optimal preconditioner performance.}

    \label{fig:pc_count}
\end{figure}

Figure~\ref{fig:pc_count} displays the distribution of preconditioners based on their optimal performance across the matrix set. As depicted, ILU-1 emerges as the most frequently optimal preconditioner, chosen in 46\% of the cases. EISENSTAT is the second most common, with 35\%. The other preconditioners were optimal less frequently, with none surpassing 14\% in frequency.  These figures suggest a dominance of ILU-1 and EISENSTAT in our data set, with other preconditioners playing a more specialized role.

\subsection{Machine Learning Models}\label{sec:ml_models}
We defined the preconditioner selection problem as a multi-label classification task, where the input is the matrix's features, and the output is a set of optimal preconditioners. We considered two types of input data: scalar features and image data. We computed the scalar and image-based features described in Section~\ref{sec:data_set}. We used the following libraries to implement the ML models: scikit-learn, TensorFlow, and Keras \cite{scikit-learn,chollet2015keras,tensorflow2015-whitepaper}.

We applied the following classifiers to the scalar feature dataset: Logistic Regression, Support Vector Classifier (SVC), Decision Tree, Random Forest, Gradient Boosting, K-Nearest Neighbors, Multi-layer Perceptron (MLP), and Gaussian Naive Bayes. While describing these traditional ML methods is beyond the scope of this article, interested readers are encouraged to consult  \cite{kotsiantis2007supervised} for a comprehensive explanation of them. We did not adjust the various parameters; therefore, each classifier used its \textit{default} parameters.

\begin{figure}
    \centering
    \includegraphics[width=0.95\linewidth]{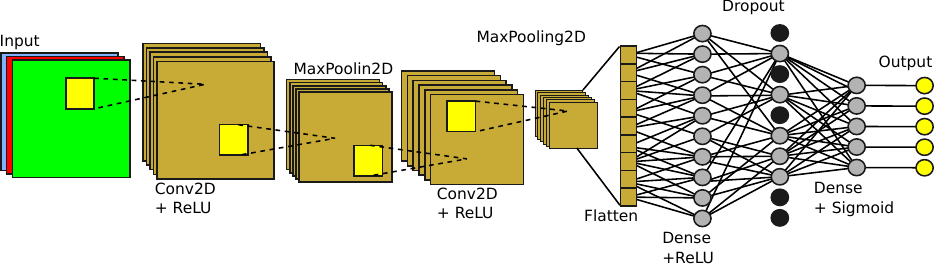}
    \caption{Convolutional Neural Network Architecture.}
    \label{fig:cnn_arch}
\end{figure}

For the RGB images, we employ a CNN. This architecture learns to identify hierarchical attributes in input images, starting with elementary attributes in the initial layers, and gradually learns to recognize complex attributes in later layers~\cite{Rawat2017DeepCN}. We apply a CNN architecture that consists of the following layers (see Figure~\ref{fig:cnn_arch}):

\begin{description}
\item[1. Input Layer:] The architecture starts with an input layer, which accepts the image data and applies a $3\times 3$ convolution operation using 32 filters. We apply a Rectified Linear Unit (ReLU) activation function to introduce non-linearity into the model. The padding is `same,' which ensures that the spatial dimensions of the input and output feature maps remain identical.

\item[2. MaxPooling Layer:] The subsequent layer is a max pooling layer with a $2\times 2$ operation on the output from the input layer. This approach reduces the spatial dimensions of the feature maps and incorporates a degree of translation invariance into the model.

\item[3. Convolutional Layer:] We apply a $3\times 3$ convolution operation with 64 filters in this layer. The ReLU activation function uses `same' padding to maintain the spatial dimensions.

\item[4. MaxPooling Layer:] Another max pooling layer is implemented with a $2\times 2$ operation, reducing the feature map dimensions.

\item[5. Flatten Layer:] The resultant feature map becomes a 1D vector. This operation is crucial for the subsequent fully connected layer, which requires its input in a 1D format.

\item[6. Fully Connected Layer:] This layer consists of 128 neurons using the ReLU activation function and accepts the flattened output from the preceding layer.

\item[7. Dropout Layer:] To prevent overfitting, a dropout layer is included, which randomly disregards 50\% of the previous layer's neurons during the training process.

\item[8. Output Layer:] Lastly, the output layer employs a sigmoid activation function to produce the final model output. The number of neurons in this layer corresponds to the number of classes (preconditioners) in the data set.
\end{description}

 Figure \ref{fig:cnn_arch} presents a graphical representation of this CNN architecture.

Multi-label classification algorithms can yield an empty set of labels (preconditioners) as an output. In such instances, we adopted the alternative of selecting the most frequent preconditioner in the data set, ILU-1. This alternative would be a natural choice for a user without additional information.

\section{Results}\label{sec:results}
The computational test and train environment was an Intel(R) Xeon(R) CPU E5-2650 v2 @ 2.60GHz processor with 64GB of RAM. The operational system was Debian GNU/Linux 11. C++ was the code language, compiled with GNU g++ (Debian 10.2.1-6) using the following directives: \mbox{-std=c++17} \mbox{-O3} \mbox{-lrt} \mbox{-fPIC} -m64 -march=native -mtune=native -fopenmp-simd. The compiler for the PETSc library was mpicxx over GNU g++ (Debian 3.3.1-5). However, the tests ran sequentially without using PETSc's MPI parallelism support, just with a single thread and one MPI process. TensorFlow (version 2.12.0)  and scikit-learn (version 1.2.2) packages \cite{scikit-learn,tensorflow2015-whitepaper} provided the machine learning environment.

We combine the scalar features and image data to create three distinct data sets. The first data set comprises the scalar features listed in Table \ref{tab:scalar_features}. The second data set includes the images. The third data set, the extended data set, contains the scalar features and flattened images (as detailed in Section~\ref{subsec:extended}).
We used the same machine-learning algorithms 
to train the scalar and the extended scalar data set.
Moreover, in the latter set, we employed Principal Component Analysis (PCA) to reduce image data dimensionality, accounting for 99\% of the variance \cite{guyon2006introduction}.

We evaluated the performance of the machine learning models using \emph{accuracy} and \emph{slowdown} as metrics. The accuracy is the proportion of correctly predicted optimal preconditioners. Specifically, let $Y_{i}$ represent the set of optimal preconditioners for matrix $i$, and $\hat{Y}_{ij}$ represent the set of preconditioners predicted by model~$j$ for matrix~$i$, then the accuracy of the model~$j$ for matrix~$i$ is given by:
\begin{equation}\label{eq:accuracy}
    \text{accuracy}_{ij} = \frac{|Y_i \cap \hat{Y}_{ij}|}{|\hat{Y}_{ij}|}.
\end{equation}

The slowdown is the ratio of the minimal execution time of the preconditioners predicted by the model to the execution time of the fastest optimal preconditioner. Let $t_{ijk}$ represent the execution time of the preconditioner $k$ predicted by model $j$ for matrix $i$, and $t_{i}^*$ represent the execution time of the fastest optimal preconditioner, then the slowdown of model $j$ for matrix $i$ is:
\begin{equation}\label{eq:slowdown}
    \text{slowdown}_{ij} = \frac{\min_k t_{ijk}}{t_{i}^*}.
\end{equation}

The metrics \emph{slowdown} and \emph{accuracy} are complementary; on the one hand, the \emph{slowdown} metric is optimistic as it considers the best of all the preconditioners predicted by the model. This aspect of \emph{slowdown} is counterbalanced by the \emph{accuracy} metric that punishes classifiers that assign too many labels (preconditioners) that are not optimal.

When evaluating the performance of each classifier, the data sets were divided into training and testing sets, maintaining a splitting of 80\%-20\%. To increase the reliability of the results and reduce possible bias from a single data splitting, we repeated the previous training-test-evaluation cycle 30 times. The intention behind this repetition was to achieve a broader perspective on the performance of the methods on different subsets of data.

\subsection{Performance on Scalar Data Set}
\begin{figure}
    \centering
    \includegraphics[width=0.95\linewidth]{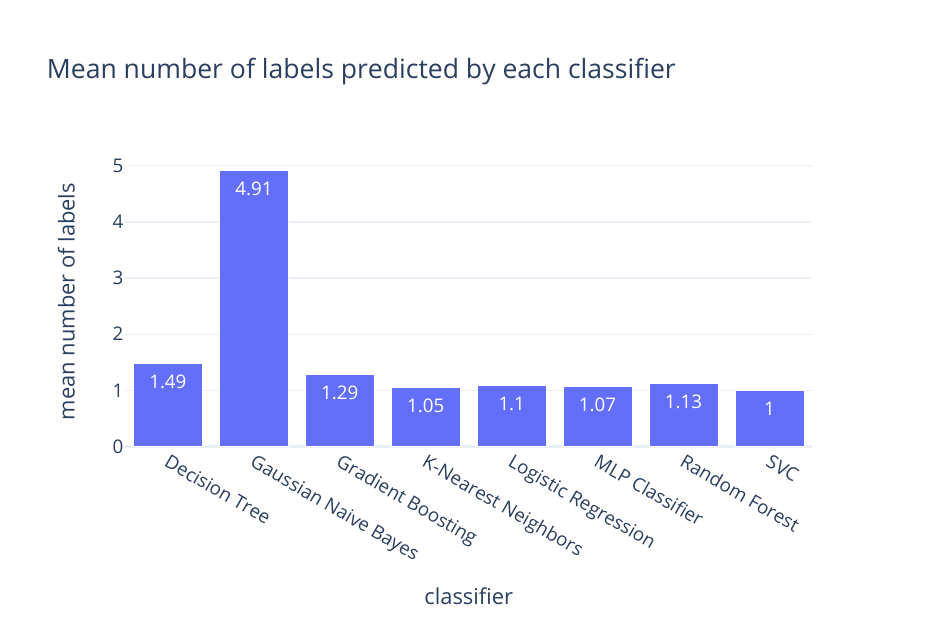}
    \caption{The mean number of labels predicted by each classifier.}
    \label{fig:MLS_num_labels}
\end{figure}

Table \ref{tab:MLS} presents the probability (relative frequency) of each classifier reaching the maximum \emph{accuracy} and an acceptable \emph{slowdown} (of 1.5 or less). Although Gaussian Naive Bayes has a high probability of 0.79 for an \emph{slowdown} below 1.5, its chances of achieving a perfect \emph{accuracy}, represented as $P(\text{accuracy} = 1) $, are remarkably low at 0.05. Figure \ref{fig:MLS_num_labels} indicates that Gaussian Naive Bayes often selects many preconditioners, with a good portion possibly not being optimal for the assigned matrices. The average number of optimal preconditioners per matrix in the data set is 1.29. In contrast, on average, the Gaussian Naive Bayes classifier selects 4.9 preconditioners with a standard deviation of 2.15. This disparity suggests that the classifier often selects non-optimal preconditioners for the given matrices; its low \emph{accuracy} of 0.05 corroborates this statement.

Conversely, Multi-layer Perceptron (MLP), Support Vector Classifier (SVC), Random Forest, and K-Nearest Neighbors exhibit balanced performance metrics. Specifically, they achieve accuracy scores exceeding 0.55 and sustain probabilities greater than 0.6 for maintaining an acceptable computational slowdown.

\begin{table}
\centering
\caption{The slowdown and accuracy of the scalar-based classifiers.}
\label{tab:MLS}
\begin{tabular}{lrr}
\toprule
classifier & P(accuracy = 1) & P(slowdown < 1.5) \\
\midrule
MLP Classifier & \textbf{0.57} & 0.69 \\
SVC & 0.56 & 0.64 \\
Random Forest & 0.55 & 0.70 \\
K-Nearest Neighbors & 0.55 & 0.66 \\
Logistic Regression & 0.48 & 0.61 \\
Gradient Boosting & 0.47 & 0.68 \\
Decision Tree & 0.39 & 0.69 \\
Gaussian Naive Bayes & 0.05 & \textbf{0.79} \\
\bottomrule
\end{tabular}
\end{table}

\begin{figure}
    \centering
    \includegraphics[width=0.95\linewidth]{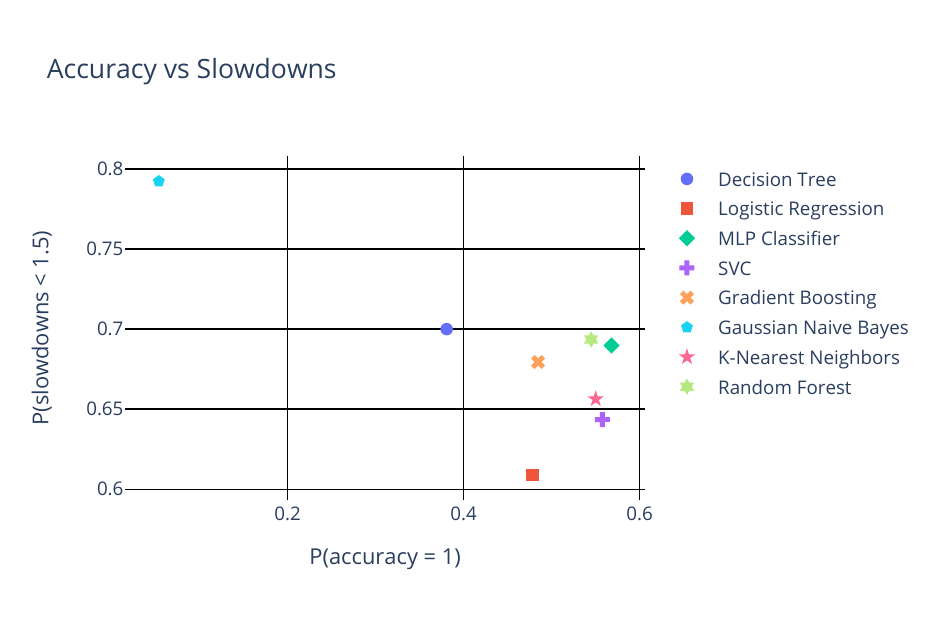}
    \caption{The performance of the scalar-based classifiers.}
    \label{fig:MLS}
\end{figure}

The plot in Figure \ref{fig:MLS} visually represents the performance trade-off between accuracy and slowdown for the scalar-based classifiers. Notably, towards the upper right, Random Forest and MLP Classifier indicate higher likelihoods of optimal accuracy and acceptable slowdown. On the other hand, Gaussian Naive Bayes is distinctly away from the main cluster, highlighting its reduced accuracy despite having a favorable slowdown probability. 

\subsection{Performance on Image Data Set}
In our study, convolutional neural networks (CNNs) were applied to an image data set. We tested four data sets comprising image resolutions of $32 \times 32$, $64 \times 64$, $128 \times 128$, and $256 \times 256$ pixels. The CNN models are identified with the suffix \_32, \_64, \_128, and \_256 to specify the associated image data set.

Table \ref{tab:CNN} presents the likelihood of maximum accuracy and slowdown of less than 1.5 for the CNN classifiers. The CNN\_256 displayed the highest probability of maximum accuracy. The CNN\_128 classifier closely trailed it with a probability of 0.88. Conversely, the CNN\_32 classifier had the lowest accuracy among the image-based models, registering a commendable 0.82 prediction probability. 
\begin{table}
\centering
\caption{Slowdown and accuracy of the image-based classifiers.}
\label{tab:CNN}
\begin{tabular}{lrr}
\toprule
classifier & \( P(\text{accuracy} = 1) \) & \( P(\text{slowdown} < 1.5) \) \\
\midrule
CNN\_256 & \textbf{0.89} & \textbf{0.94} \\
CNN\_128 & 0.88 & \textbf{0.94} \\
CNN\_64 & 0.86 & 0.93 \\
CNN\_32 & 0.82 & 0.89 \\
\bottomrule
\end{tabular}
\end{table}

Regarding slowdown metrics, all four classifiers demonstrated impressive performances, with the CNN\_256 and CNN\_128 classifiers showcasing a high probability of 0.94 for maintaining a slowdown factor below 1.5. The CNN\_64 was behind with a probability of 0.93, while the CNN\_32 classifier exhibited a probability of 0.89.

Finally, the same table shows that while a rise in image resolution can enhance results, accuracy and slowdown eventually plateau. This observation suggests that using larger images does not offer additional benefits.

\subsection{Performance on Extended Data Set}

We compiled an extended data set combining scalar features with vectors derived from flattened images to bolster our analysis and emphasize the benefits of image encoding for matrix attributes. This amalgamated data set directly compares image-based topological information and its scalar counterparts. We generated four distinct extended data sets and subjected each to the scalar classifiers. We adopted suffix notations \_32, \_64, \_128, and \_256 to denote the image resolution linked to each data set and classifier.

Due to the large number of configurations, we streamlined our focus on the top 10 setups for each metric, postponing the detailed results to the Appendix -  Mixed Models \ref{appendix:mixed_models}.

Table \ref{tab:MLM_merged_top10_accuracy} details the probability of obtaining maximum accuracy across the top 10 classifiers. Notably, the Random Forest classifier, especially when paired with the 256x256 image resolution, led the list with an accuracy rate of 0.60. Across varied image resolutions, Random Forest consistently secured its top-tier status. 

\begin{table}
\centering
\caption{Top 10 classifiers by accuracy in the extended scalar-based feature set.}
\label{tab:MLM_merged_top10_accuracy}
\begin{tabular}{lr}
\toprule
classifier & P(accuracy = 1) \\
\midrule
Random Forest\_256 & \textbf{0.60} \\
Random Forest\_32 & 0.56 \\
Random Forest\_64 & 0.56 \\
Random Forest\_128 & 0.55 \\
SVC\_256 & 0.53 \\
Gradient Boosting\_256 & 0.52 \\
SVC\_32 & 0.52 \\
K-Nearest Neighbors\_32 & 0.51 \\
SVC\_128 & 0.51 \\
SVC\_64 & 0.50 \\
\bottomrule
\end{tabular}
\end{table}

Pivoting to the slowdown metrics highlighted in Table \ref{tab:MLM_merged_top10_slowdowns}, the Random Forest classifier, particularly the 256x256 resolution variant, once again exhibited superior performance. It demonstrated a 0.71 probability of achieving a slowdown of less than 1.5. Several other classifiers followed closely, like the Gradient Boosting and SVC across different resolutions.

\begin{table}
\centering
\caption{Top 10 classifiers by slowdown in the extended scalar-based feature set.}
\label{tab:MLM_merged_top10_slowdowns}
\begin{tabular}{lr}
\toprule
classifier & P(slowdown < 1.5) \\
\midrule
Random Forest\_256 & \textbf{0.71} \\
Gradient Boosting\_256 & \textbf{0.71} \\
Random Forest\_64 & 0.68 \\
Random Forest\_32 & 0.68 \\
Random Forest\_128 & 0.68 \\
K-Nearest Neighbors\_32 & 0.65 \\
SVC\_32 & 0.61 \\
SVC\_256 & 0.61 \\
SVC\_128 & 0.61 \\
SVC\_64 & 0.59 \\
\bottomrule
\end{tabular}
\end{table}

\begin{figure}
    \centering
    \includegraphics[width=0.95\linewidth]{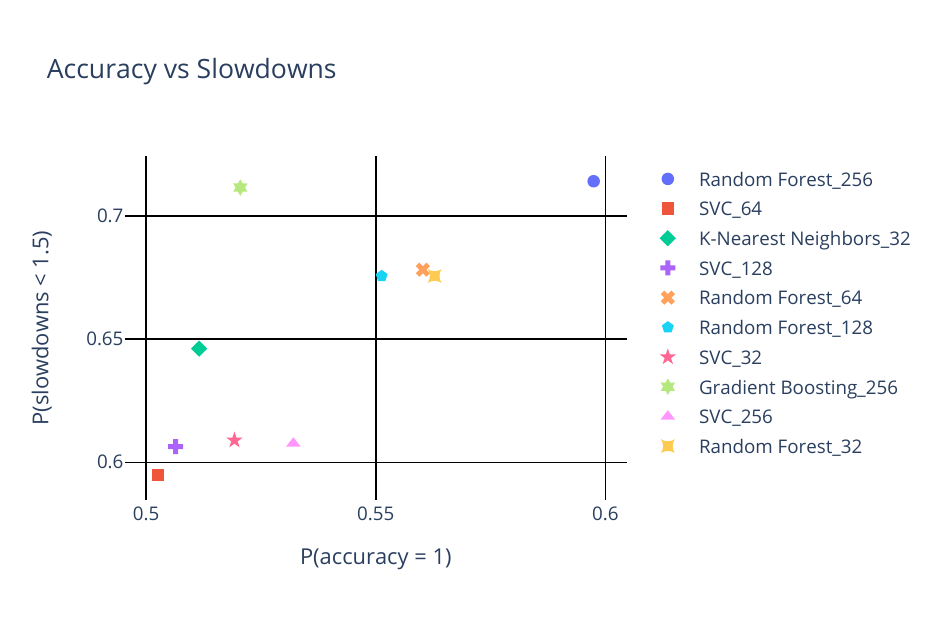}
    \caption{The extended classifiers' performance is based on the probability of achieving perfect accuracy and the likelihood of obtaining a slowdown of less than 1.5.}
    \label{fig:MLM_merged}
\end{figure}

Figure \ref{fig:MLM_merged} provides a comparative visualization of the performance across the extended models. The plot demonstrates the Random Forest\_256 model's dominance, positioning it as the standout choice in the evaluated set.

\subsection{Comparison of Scalar, Image, and Extended Models}
In this final result section, we perform a comparative analysis to identify the best classifier configurations from the three previously discussed sections. The objective is to highlight the most effective approach for achieving high accuracy and acceptable slowdown. As a reference point, we included a Benchmark classifier, which selects the most frequently occurring preconditioner in the data set. In our case, this preconditioner was ILU-1.

Table \ref{tab:ALL_merged} presents the best models' slowdown and accuracy across each feature set. The CNN configurations consistently exhibit exceptional performance. Specifically, the CNN\_256 model achieved the highest accuracy of 0.89 and a remarkable slowdown metric of 0.94, making it the standout performer among all the models. On the non-CNN models, Random Forest\_256 recorded an accuracy of 0.60 and a slowdown metric of 0.71, closely followed by the MLP Classifier, which achieved an accuracy of 0.57 and a slowdown metric of 0.69. The data suggests the efficacy of image representations of matrix sparsity in advancing preconditioner selection algorithms.

\begin{table}
    \centering
    \caption{The slowdown and accuracy of the best models of each data set of features.}
    \label{tab:ALL_merged}
    \begin{tabular}{lrr}
    \toprule
    classifier & P(accuracy = 1) & P(slowdown < 1.5) \\
    \midrule
    Benchmark (ILU-1) & 0.49 & 0.58 \\
    MLP Classifier & 0.57 & 0.69 \\
    SVC & 0.56 & 0.64 \\
    Random Forest & 0.55 & 0.70 \\
    CNN\_256 & \textbf{0.89}& \textbf{0.94}\\
    CNN\_128 & 0.88 & \textbf{0.94}\\
    CNN\_64 & 0.86 & 0.93 \\
    CNN\_32 & 0.82 & 0.89 \\
    Random Forest\_256 & 0.60 & 0.71 \\
    Gradient Boosting\_256 & 0.52 & 0.71 \\
    \bottomrule
    \end{tabular}
\end{table}

\begin{figure}
    \centering
    \includegraphics[width=0.95\linewidth]{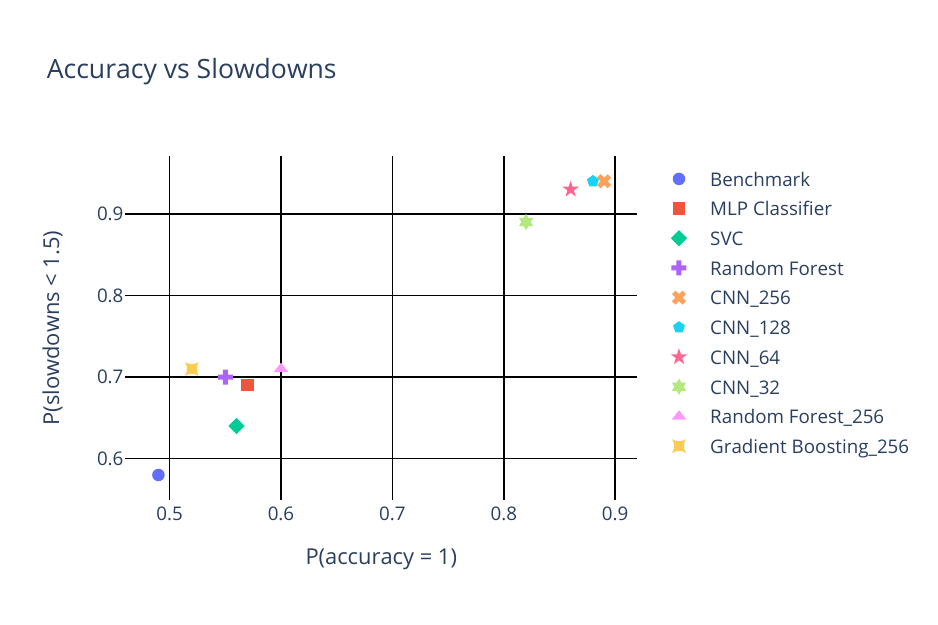}
    \caption{The performance of the best models of each data set of features.}
    \label{fig:ALL_merged}
\end{figure}

Building upon Table \ref{tab:ALL_merged}, Figure \ref{fig:ALL_merged} visually maps each classifier's relationship between accuracy and slowdown. The CNN configurations are clustered in the top-right quadrant, meaning superior accuracy and low slowdown performance. This graphical representation reiterates the dominance of CNN-based models over others, with CNN\_256 being the standout performer. Conversely, the Benchmark model is positioned in the bottom-left, underlining its relatively inferior performance. The non-CNN models, including the MLP Classifier, Random Forest, and SVC, are scattered around the central region, further emphasizing the distinction between the effectiveness of the two approaches.

\section{Conclusions}\label{sec:conclusions}

Effective solvers tailored for large, sparse linear systems are crucial in high-performance computing. These solvers primarily use iterative methods, and their efficacy largely hinges on using a preconditioner. However, choosing the suitable preconditioner is intricate since there is no one-size-fits-all approach for all matrices and applications. Choosing the suitable preconditioner depends on several factors, with matrix sparsity being particularly important. Although bandwidth and the number of non-zero elements provide some insights, they cannot capture the full complexity of sparsity patterns.

Our research extends the Yamada et al. \cite{yamada2018preconditioner} work; we also encoded matrix features as RGB images. We employed a convolutional neural network (CNN) to streamline the preconditioner selection, framing it as a multi-class classification task. We undertook a comparative analysis to highlight the merits of this image-based approach. The CNN model trained using images was juxtaposed against models based on scalar attributes, encompassing density, matrix size, eigenvalues, and conditioning. Additionally, an examination was conducted on 
an extended scalar data set,
amalgamating scalar attributes with flattened vectors of images.

Our empirical investigations, harnessing 126 SPD matrices extracted from the SuiteSparse Matrix Collection, yielded compelling insights. The proposed CNN model registered a significant enhancement, clocking in a 32\% surge in the likelihood of maximum accuracy and a 25\% improvement in the probability of acceptable slowdown. Such findings reaffirm the potential of encoding matrix as an image, positioning it as a prospective alternative in preconditioner selection.

Charting the path ahead, we envision leveraging CNN to fine-tune other parameters of the iterative methods like the specific solver type, tolerance parameters, and so forth, paving the way for more optimized and informed computational solutions in HPC applications. 

\section{Acknowledgements}
This work was supported by the Petróleo Brasileiro S.A. (Petrobras) under grant PT-200.20.00125.

\begin{landscape}
\section*{Appendix - Matrix Data Set}\label{appendix:matrix_data_set}
\begin{longtable}{lccccccccc}
  \toprule
  Mat & NRows & NNZ & Density & RowNNZ & MaxEigs & MinEigs & Condest & DDom & DDeg \\
  \midrule
  \endfirsthead
  \multicolumn{10}{c}{\tablename\ \thetable{} -- continued from previous page} \\
   \midrule
  Mat & NRows & NNZ & Density & RowNNZ & MaxEigs & MinEigs & Condest & DDom & DDeg \\
  \midrule
  \endhead
   \midrule
  \multicolumn{10}{r}{Continued on next page...} \\
  \endfoot
   \midrule
  \endlastfoot
denormal & 8.94E+04 & 1.16E+06 & 1.45E-04 & 1.29E+01 & 1.28E-01 & 2.77E-19 & 1.07E+17 & 0.00E+00 & 4.29E-01 \\
2cubes\_sphere & 1.01E+05 & 1.65E+06 & 1.60E-04 & 1.62E+01 & 2.83E+10 & 3.22E+01 & 2.94E+09 & 7.13E-01 & 3.33E-01 \\
gyro & 1.74E+04 & 1.02E+06 & 3.39E-03 & 5.88E+01 & 3.66E+09 & 3.34E+00 & 4.44E+09 & 2.21E-02 & 8.29E-02 \\
bcsstk24 & 3.56E+03 & 1.60E+05 & 1.26E-02 & 4.49E+01 & 3.07E+13 & 1.57E+02 & 6.37E+11 & 3.23E-01 & 4.40E-03 \\
bcsstk26 & 1.92E+03 & 3.03E+04 & 8.21E-03 & 1.58E+01 & 1.58E+11 & 9.54E+02 & 4.39E+08 & 0.00E+00 & 5.08E-04 \\
bcsstk16 & 4.88E+03 & 2.90E+05 & 1.22E-02 & 5.95E+01 & 4.94E+09 & 1.00E+00 & 7.01E+09 & 1.58E-02 & 2.67E-01 \\
msc01440 & 1.44E+03 & 4.50E+04 & 2.17E-02 & 3.12E+01 & 1.36E+06 & 4.10E-01 & 7.00E+06 & 0.00E+00 & 6.88E-04 \\
fv1 & 9.60E+03 & 8.53E+04 & 9.24E-04 & 8.88E+00 & 6.54E+00 & 5.12E-01 & 1.28E+01 & 0.00E+00 & 1.04E-04 \\
plbuckle & 1.28E+03 & 3.06E+04 & 1.86E-02 & 2.39E+01 & 2.10E+06 & 1.63E+00 & 4.13E+06 & 0.00E+00 & 7.80E-04 \\
bcsstk21 & 3.60E+03 & 2.66E+04 & 2.05E-03 & 7.39E+00 & 1.27E+08 & 7.21E+00 & 4.50E+07 & 0.00E+00 & 2.78E-04 \\
Dubcova3 & 1.47E+05 & 3.64E+06 & 1.69E-04 & 2.48E+01 & 4.80E+00 & 1.20E-03 & 1.15E+04 & 0.00E+00 & 3.94E-01 \\
bcsstk09 & 1.08E+03 & 1.84E+04 & 1.57E-02 & 1.70E+01 & 6.76E+07 & 7.10E+03 & 3.10E+04 & 0.00E+00 & 9.23E-04 \\
Trefethen\_150 & 1.50E+02 & 2.04E+03 & 9.07E-02 & 1.36E+01 & 8.64E+02 & 1.12E+00 & 1.13E+03 & 0.00E+00 & 6.71E-03 \\
Trefethen\_300 & 3.00E+02 & 4.68E+03 & 5.20E-02 & 1.56E+01 & 1.99E+03 & 1.12E+00 & 2.58E+03 & 0.00E+00 & 3.34E-03 \\
bcsstk18 & 1.19E+04 & 1.49E+05 & 1.04E-03 & 1.25E+01 & 4.30E+10 & 1.24E-01 & 6.49E+11 & 4.21E-01 & 4.51E-03 \\
LFAT5000 & 2.00E+04 & 8.00E+04 & 2.00E-04 & 4.00E+00 & 7.26E+13 & 1.50E-04 & 4.83E+17 & 0.00E+00 & 1.05E-05 \\
nos2 & 9.57E+02 & 4.14E+03 & 4.52E-03 & 4.32E+00 & 1.57E+11 & 3.08E+01 & 6.34E+09 & 0.00E+00 & 9.99E-04 \\
Kuu & 7.10E+03 & 3.40E+05 & 6.74E-03 & 4.79E+01 & 5.41E+01 & 3.43E-03 & 3.26E+04 & 0.00E+00 & 1.63E-01 \\
cvxbqp1 & 5.00E+04 & 3.50E+05 & 1.40E-04 & 7.00E+00 & 4.83E+05 & 3.95E-02 & 1.22E+07 & 4.00E-05 & 5.00E-01 \\
crystm01 & 4.88E+03 & 1.05E+05 & 4.43E-03 & 2.16E+01 & 5.31E-12 & 2.32E-14 & 4.21E+02 & 0.00E+00 & 2.05E-04 \\
bundle1 & 1.06E+04 & 7.71E+05 & 6.88E-03 & 7.28E+01 & 6.43E+12 & 4.83E+08 & 1.33E+04 & 3.07E-01 & 2.58E-02 \\
shallow\_water1 & 8.19E+04 & 3.28E+05 & 4.88E-05 & 4.00E+00 & 2.03E+10 & 5.79E+09 & 3.63E+00 & 1.00E+00 & 1.83E+00 \\
Chem97ZtZ & 2.54E+03 & 7.36E+03 & 1.14E-03 & 2.90E+00 & 1.33E+03 & 5.39E+00 & 4.63E+02 & 0.00E+00 & 3.94E-04 \\
bcsstk06 & 4.20E+02 & 7.86E+03 & 4.46E-02 & 1.87E+01 & 3.49E+09 & 4.61E+02 & 1.22E+07 & 0.00E+00 & 2.14E-03 \\
obstclae & 4.00E+04 & 1.98E+05 & 1.24E-04 & 4.94E+00 & 8.20E+00 & 2.00E-01 & 4.10E+01 & 1.00E+00 & 1.05E+00 \\
cfd1 & 7.07E+04 & 1.83E+06 & 3.66E-04 & 2.58E+01 & 6.80E+00 & 2.00E-05 & 1.34E+06 & 0.00E+00 & 1.26E-01 \\
bcsstk25 & 1.54E+04 & 2.52E+05 & 1.06E-03 & 1.63E+01 & 1.06E+15 & 2.40E+02 & 1.28E+13 & 3.42E-01 & 2.79E-03 \\
cant & 6.25E+04 & 4.01E+06 & 1.03E-03 & 6.42E+01 & 1.94E+04 & 7.46E-07 & 5.76E+10 & 3.89E-03 & 5.91E-02 \\
fv2 & 9.80E+03 & 8.70E+04 & 9.06E-04 & 8.88E+00 & 6.54E+00 & 5.12E-01 & 1.28E+01 & 0.00E+00 & 1.02E-04 \\
gyro\_m & 1.74E+04 & 3.40E+05 & 1.13E-03 & 1.96E+01 & 1.89E-07 & 7.54E-14 & 1.12E+07 & 0.00E+00 & 1.67E-01 \\
apache1 & 8.08E+04 & 5.42E+05 & 8.30E-05 & 6.71E+00 & 3.63E+05 & 9.09E-02 & 3.99E+06 & 2.70E-01 & 1.00E+00 \\
minsurfo & 4.08E+04 & 2.04E+05 & 1.22E-04 & 4.99E+00 & 8.10E+00 & 1.00E-01 & 8.11E+01 & 1.00E+00 & 1.00E+00 \\
wathen120 & 3.64E+04 & 5.66E+05 & 4.26E-04 & 1.55E+01 & 3.69E+02 & 1.43E-01 & 4.05E+03 & 0.00E+00 & 1.71E-01 \\
cbuckle & 1.37E+04 & 6.77E+05 & 3.61E-03 & 4.94E+01 & 7.22E+04 & 2.19E-03 & 8.09E+07 & 4.53E-02 & 1.04E-02 \\
bcsstk15 & 3.95E+03 & 1.18E+05 & 7.56E-03 & 2.98E+01 & 6.54E+09 & 1.00E+00 & 7.97E+09 & 0.00E+00 & 2.50E-04 \\
ex10hs & 2.55E+03 & 5.73E+04 & 8.83E-03 & 2.25E+01 & 5.27E+07 & 9.61E-05 & 1.27E+12 & 0.00E+00 & 3.85E-04 \\
pdb1HYS & 3.64E+04 & 4.34E+06 & 3.28E-03 & 1.19E+02 & 3.52E+02 & 9.97E-10 & 2.46E+12 & 0.00E+00 & 1.05E-01 \\
LFAT5000\_E & 2.00E+04 & 8.00E+04 & 2.00E-04 & 4.00E+00 & 7.26E+11 & 1.50E-06 & 4.83E+17 & 0.00E+00 & 1.05E-05 \\
ecology2 & 1.00E+06 & 5.00E+06 & 5.00E-06 & 5.00E+00 & 8.17E+01 & 1.23E-06 & 6.66E+07 & 7.00E-06 & 1.00E+00 \\
bcsstk07 & 4.20E+02 & 7.86E+03 & 4.46E-02 & 1.87E+01 & 3.49E+09 & 4.61E+02 & 1.22E+07 & 0.00E+00 & 2.14E-03 \\
t2dah\_e & 1.14E+04 & 1.76E+05 & 1.34E-03 & 1.54E+01 & 2.29E-05 & 3.16E-14 & 1.37E+09 & 0.00E+00 & 1.00E-01 \\
ted\_B\_unscaled & 1.06E+04 & 1.45E+05 & 1.29E-03 & 1.36E+01 & 1.00E+00 & 7.82E-12 & 2.04E+11 & 6.29E-01 & 2.59E-01 \\
s2rmq4m1 & 5.49E+03 & 2.63E+05 & 8.74E-03 & 4.80E+01 & 6.87E+04 & 3.87E-04 & 3.22E+08 & 3.64E-04 & 4.15E-03 \\
Trefethen\_2000 & 2.00E+03 & 4.19E+04 & 1.05E-02 & 2.10E+01 & 1.74E+04 & 1.12E+00 & 2.25E+04 & 0.00E+00 & 5.00E-04 \\
jnlbrng1 & 4.00E+04 & 1.99E+05 & 1.24E-04 & 4.98E+00 & 1.84E+01 & 1.01E-01 & 1.87E+02 & 1.00E+00 & 1.01E+00 \\
bcsstk14 & 1.81E+03 & 6.35E+04 & 1.95E-02 & 3.51E+01 & 1.19E+10 & 1.00E+00 & 1.31E+10 & 0.00E+00 & 5.08E-04 \\
qa8fm & 6.61E+04 & 1.66E+06 & 3.80E-04 & 2.51E+01 & 9.46E-04 & 1.30E-05 & 1.10E+02 & 0.00E+00 & 4.21E-01 \\
gyro\_k & 1.74E+04 & 1.02E+06 & 3.39E-03 & 5.88E+01 & 3.66E+09 & 3.34E+00 & 4.44E+09 & 2.21E-02 & 8.29E-02 \\
thermomech\_TC & 1.02E+05 & 7.12E+05 & 6.82E-05 & 6.97E+00 & 3.05E-02 & 4.54E-04 & 1.24E+02 & 5.02E-01 & 1.00E+00 \\
ex10 & 2.41E+03 & 5.48E+04 & 9.44E-03 & 2.28E+01 & 6.96E+07 & 7.65E-05 & 2.06E+12 & 0.00E+00 & 4.04E-04 \\
gr\_30\_30 & 9.00E+02 & 7.74E+03 & 9.56E-03 & 8.60E+00 & 1.20E+01 & 6.15E-02 & 3.77E+02 & 0.00E+00 & 1.11E-03 \\
Dubcova2 & 6.50E+04 & 1.03E+06 & 2.44E-04 & 1.58E+01 & 4.80E+00 & 1.20E-03 & 1.04E+04 & 1.23E-04 & 4.49E-01 \\
gridgena & 4.90E+04 & 5.12E+05 & 2.14E-04 & 1.05E+01 & 2.77E+04 & 1.78E-01 & 6.26E+05 & 9.24E-01 & 4.82E-01 \\
mhdb416 & 4.16E+02 & 2.31E+03 & 1.34E-02 & 5.56E+00 & 2.20E+00 & 5.50E-10 & 5.05E+09 & 0.00E+00 & 1.55E-03 \\
Trefethen\_20000b & 2.00E+04 & 5.54E+05 & 1.39E-03 & 2.77E+01 & 2.25E+05 & 2.34E+00 & 1.33E+05 & 1.00E+00 & 2.00E-01 \\
nos1 & 2.37E+02 & 1.02E+03 & 1.81E-02 & 4.29E+00 & 2.46E+09 & 1.23E+02 & 2.53E+07 & 0.00E+00 & 4.07E-03 \\
torsion1 & 4.00E+04 & 1.98E+05 & 1.24E-04 & 4.94E+00 & 8.20E+00 & 2.00E-01 & 4.10E+01 & 1.00E+00 & 1.05E+00 \\
bcsstm12 & 1.47E+03 & 1.97E+04 & 9.06E-03 & 1.33E+01 & 1.34E+01 & 2.12E-05 & 8.88E+05 & 0.00E+00 & 6.76E-04 \\
1138\_bus & 1.14E+03 & 4.05E+03 & 3.13E-03 & 3.56E+00 & 3.01E+04 & 3.52E-03 & 1.23E+07 & 0.00E+00 & 8.79E-04 \\
msc04515 & 4.52E+03 & 9.77E+04 & 4.79E-03 & 2.16E+01 & 3.15E+10 & 1.39E+04 & 5.48E+06 & 0.00E+00 & 2.21E-04 \\
bcsstk10 & 1.09E+03 & 2.21E+04 & 1.87E-02 & 2.03E+01 & 4.47E+07 & 8.54E+01 & 1.32E+06 & 0.00E+00 & 9.20E-04 \\
nd3k & 9.00E+03 & 3.28E+06 & 4.05E-02 & 3.64E+02 & 1.28E+02 & 8.14E-06 & 5.68E+07 & 0.00E+00 & 7.03E-02 \\
apache2 & 7.15E+05 & 4.82E+06 & 9.42E-06 & 6.74E+00 & 1.59E+05 & 2.99E-02 & 5.32E+06 & 3.83E-01 & 8.96E-01 \\
mesh2em5 & 3.06E+02 & 2.02E+03 & 2.16E-02 & 6.59E+00 & 2.47E+02 & 1.00E+00 & 2.93E+02 & 0.00E+00 & 3.28E-03 \\
mesh3em5 & 2.89E+02 & 1.38E+03 & 1.65E-02 & 4.76E+00 & 4.97E+00 & 1.00E+00 & 5.00E+00 & 0.00E+00 & 3.48E-03 \\
662\_bus & 6.62E+02 & 2.47E+03 & 5.65E-03 & 3.74E+00 & 4.01E+03 & 5.05E-03 & 8.27E+05 & 0.00E+00 & 1.51E-03 \\
cfd2 & 1.23E+05 & 3.09E+06 & 2.02E-04 & 2.50E+01 & 7.69E+00 & 5.06E-06 & 3.73E+06 & 0.00E+00 & 1.28E-01 \\
shallow\_water2 & 8.19E+04 & 3.28E+05 & 4.88E-05 & 4.00E+00 & 6.10E+10 & 5.80E+09 & 1.13E+01 & 1.00E+00 & 1.21E+00 \\
685\_bus & 6.85E+02 & 3.25E+03 & 6.92E-03 & 4.74E+00 & 2.62E+04 & 6.19E-02 & 5.31E+05 & 0.00E+00 & 1.46E-03 \\
fv3 & 9.80E+03 & 8.70E+04 & 9.06E-04 & 8.88E+00 & 8.72E+00 & 1.97E-03 & 4.42E+03 & 0.00E+00 & 1.02E-04 \\
Andrews & 6.00E+04 & 7.60E+05 & 2.11E-04 & 1.27E+01 & 3.65E+01 & 3.01E-16 & 2.80E+17 & 0.00E+00 & 1.00E+00 \\
s2rmt3m1 & 5.49E+03 & 2.18E+05 & 7.22E-03 & 3.97E+01 & 9.67E+04 & 3.87E-04 & 4.84E+08 & 3.64E-04 & 3.77E-03 \\
s1rmt3m1 & 5.49E+03 & 2.18E+05 & 7.22E-03 & 3.97E+01 & 9.67E+05 & 3.80E-01 & 5.35E+06 & 3.64E-04 & 9.29E-03 \\
LFAT5000\_M & 2.00E+04 & 8.00E+04 & 2.00E-04 & 4.00E+00 & 1.14E-04 & 8.86E-16 & 1.29E+11 & 0.00E+00 & 1.74E-05 \\
msc00726 & 7.26E+02 & 3.45E+04 & 6.55E-02 & 4.75E+01 & 4.22E+09 & 1.01E+04 & 1.17E+06 & 0.00E+00 & 1.38E-03 \\
sts4098 & 4.10E+03 & 7.24E+04 & 4.31E-03 & 1.77E+01 & 3.07E+08 & 1.41E+00 & 4.44E+08 & 0.00E+00 & 2.38E-04 \\
nos6 & 6.75E+02 & 3.26E+03 & 7.14E-03 & 4.82E+00 & 7.65E+06 & 9.56E-01 & 8.00E+06 & 0.00E+00 & 1.48E-03 \\
bcsstk22 & 1.38E+02 & 6.96E+02 & 3.65E-02 & 5.04E+00 & 5.85E+06 & 5.28E+01 & 1.67E+05 & 0.00E+00 & 6.94E-03 \\
bcsstk34 & 5.88E+02 & 2.14E+04 & 6.19E-02 & 3.64E+01 & 3.98E+07 & 1.40E+03 & 5.09E+04 & 0.00E+00 & 1.70E-03 \\
Trefethen\_500 & 5.00E+02 & 8.48E+03 & 3.39E-02 & 1.70E+01 & 3.57E+03 & 1.12E+00 & 4.63E+03 & 0.00E+00 & 2.00E-03 \\
494\_bus & 4.94E+02 & 1.67E+03 & 6.83E-03 & 3.37E+00 & 3.00E+04 & 1.24E-02 & 3.89E+06 & 0.00E+00 & 2.02E-03 \\
bcsstk27 & 1.22E+03 & 5.61E+04 & 3.75E-02 & 4.59E+01 & 3.47E+06 & 1.44E+02 & 7.71E+04 & 0.00E+00 & 8.01E-04 \\
bodyy6 & 1.94E+04 & 1.34E+05 & 3.58E-04 & 6.93E+00 & 7.81E+04 & 1.02E+00 & 9.80E+04 & 1.00E+00 & 1.00E+00 \\
bcsstk12 & 1.47E+03 & 3.42E+04 & 1.58E-02 & 2.32E+01 & 6.56E+08 & 2.96E+00 & 5.25E+08 & 0.00E+00 & 6.75E-04 \\
parabolic\_fem & 5.26E+05 & 3.67E+06 & 1.33E-05 & 6.99E+00 & 8.03E-01 & 3.80E-06 & 2.11E+05 & 1.00E+00 & 1.00E+00 \\
thermal1 & 8.27E+04 & 5.74E+05 & 8.41E-05 & 6.95E+00 & 7.59E+00 & 2.39E-05 & 4.96E+05 & 2.03E-01 & 7.55E-01 \\
nos7 & 7.29E+02 & 4.62E+03 & 8.69E-03 & 6.33E+00 & 9.86E+06 & 4.15E-03 & 4.05E+09 & 0.00E+00 & 1.37E-03 \\
Trefethen\_20000 & 2.00E+04 & 5.54E+05 & 1.39E-03 & 2.77E+01 & 3.26E+05 & 1.12E+00 & 2.91E+05 & 1.00E+00 & 1.33E-01 \\
ted\_B & 1.06E+04 & 1.45E+05 & 1.29E-03 & 1.36E+01 & 1.00E+00 & 5.29E-08 & 3.02E+07 & 6.29E-01 & 2.59E-01 \\
bcsstk03 & 1.12E+02 & 6.40E+02 & 5.10E-02 & 5.71E+00 & 2.00E+11 & 2.94E+04 & 9.41E+06 & 0.00E+00 & 5.30E-03 \\
gyro\_M & 1.74E+04 & 3.40E+05 & 1.13E-03 & 1.96E+01 & 1.89E-07 & 7.54E-14 & 1.12E+07 & 0.00E+00 & 1.67E-01 \\
nos5 & 4.68E+02 & 5.17E+03 & 2.36E-02 & 1.11E+01 & 5.82E+05 & 5.29E+01 & 2.92E+04 & 0.00E+00 & 2.13E-03 \\
boneS01\_M & 1.27E+05 & 2.24E+06 & 1.38E-04 & 1.76E+01 & 1.79E-06 & 1.42E-08 & 1.69E+02 & 0.00E+00 & 4.21E-01 \\
bcsstk13 & 2.00E+03 & 8.39E+04 & 2.09E-02 & 4.19E+01 & 3.11E+12 & 2.84E+02 & 4.57E+10 & 0.00E+00 & 4.77E-04 \\
mesh2e1 & 3.06E+02 & 2.02E+03 & 2.16E-02 & 6.59E+00 & 3.81E+02 & 1.31E+00 & 4.17E+02 & 0.00E+00 & 3.28E-03 \\
nasa2146 & 2.15E+03 & 7.22E+04 & 1.57E-02 & 3.37E+01 & 3.27E+07 & 1.90E+04 & 3.44E+03 & 0.00E+00 & 4.65E-04 \\
bcsstk19 & 8.17E+02 & 6.85E+03 & 1.03E-02 & 8.39E+00 & 1.92E+14 & 1.43E+03 & 2.81E+11 & 0.00E+00 & 9.53E-04 \\
nos3 & 9.60E+02 & 1.58E+04 & 1.72E-02 & 1.65E+01 & 6.90E+02 & 1.83E-02 & 7.35E+04 & 0.00E+00 & 1.04E-03 \\
nasasrb & 5.49E+04 & 2.68E+06 & 8.89E-04 & 4.88E+01 & 2.65E+09 & 4.74E+00 & 1.48E+09 & 5.18E-03 & 6.30E-03 \\
bcsstk28 & 4.41E+03 & 2.19E+05 & 1.13E-02 & 4.97E+01 & 7.70E+08 & 8.14E-01 & 6.28E+09 & 0.00E+00 & 2.98E-02 \\
boneS01 & 1.27E+05 & 5.52E+06 & 3.41E-04 & 4.34E+01 & 4.85E+04 & 2.85E-03 & 4.22E+07 & 1.57E-05 & 2.29E-01 \\
Trefethen\_200b & 1.99E+02 & 2.87E+03 & 7.25E-02 & 1.44E+01 & 1.22E+03 & 2.34E+00 & 7.27E+02 & 0.00E+00 & 5.05E-03 \\
s1rmq4m1 & 5.49E+03 & 2.62E+05 & 8.71E-03 & 4.78E+01 & 6.87E+05 & 3.80E-01 & 3.58E+06 & 3.64E-04 & 6.46E-03 \\
thermomech\_dM & 2.04E+05 & 1.42E+06 & 3.41E-05 & 6.97E+00 & 6.11E-05 & 9.08E-07 & 1.25E+02 & 5.02E-01 & 1.00E+00 \\
bodyy4 & 1.75E+04 & 1.22E+05 & 3.95E-04 & 6.93E+00 & 8.32E+02 & 1.03E+00 & 1.02E+03 & 0.00E+00 & 5.70E-05 \\
Trefethen\_200 & 2.00E+02 & 2.89E+03 & 7.22E-02 & 1.44E+01 & 1.22E+03 & 1.12E+00 & 1.59E+03 & 0.00E+00 & 5.03E-03 \\
olafu & 1.61E+04 & 1.02E+06 & 3.89E-03 & 6.29E+01 & 9.48E+11 & 1.25E+00 & 2.25E+12 & 2.85E-03 & 5.51E-03 \\
bcsstk08 & 1.07E+03 & 1.30E+04 & 1.12E-02 & 1.21E+01 & 7.66E+10 & 2.95E+03 & 4.73E+07 & 0.00E+00 & 8.34E-04 \\
bcsstm07 & 4.20E+02 & 7.25E+03 & 4.11E-02 & 1.73E+01 & 2.51E+03 & 3.30E-01 & 1.34E+04 & 0.00E+00 & 2.22E-03 \\
G2\_circuit & 1.50E+05 & 7.27E+05 & 3.23E-05 & 4.84E+00 & 2.27E+04 & 1.79E-03 & 1.98E+07 & 4.97E-01 & 1.00E+00 \\
mesh3e1 & 2.89E+02 & 1.38E+03 & 1.65E-02 & 4.76E+00 & 8.93E+00 & 1.00E+00 & 9.00E+00 & 0.00E+00 & 3.48E-03 \\
Dubcova1 & 1.61E+04 & 2.53E+05 & 9.73E-04 & 1.57E+01 & 4.80E+00 & 4.81E-03 & 2.62E+03 & 4.96E-04 & 4.49E-01 \\
wathen100 & 3.04E+04 & 4.72E+05 & 5.10E-04 & 1.55E+01 & 3.70E+02 & 6.36E-02 & 8.25E+03 & 0.00E+00 & 1.71E-01 \\
Pres\_Poisson & 1.48E+04 & 7.16E+05 & 3.26E-03 & 4.83E+01 & 2.60E+01 & 1.28E-05 & 3.20E+06 & 0.00E+00 & 1.92E-01 \\
bodyy5 & 1.86E+04 & 1.29E+05 & 3.73E-04 & 6.93E+00 & 8.04E+03 & 1.02E+00 & 9.98E+03 & 0.00E+00 & 5.38E-05 \\
crystm03 & 2.47E+04 & 5.84E+05 & 9.57E-04 & 2.36E+01 & 9.79E-13 & 3.71E-15 & 4.68E+02 & 0.00E+00 & 4.03E-01 \\
aft01 & 8.20E+03 & 1.26E+05 & 1.87E-03 & 1.53E+01 & 1.00E+15 & 1.53E-04 & 9.36E+18 & 0.00E+00 & 1.22E-04 \\
smt & 2.57E+04 & 3.75E+06 & 5.67E-03 & 1.46E+02 & 5.30E+06 & 3.34E-03 & 6.13E+09 & 0.00E+00 & 6.60E-02 \\
nos4 & 1.00E+02 & 5.94E+02 & 5.94E-02 & 5.94E+00 & 8.49E-01 & 5.38E-04 & 2.70E+03 & 0.00E+00 & 1.00E-02 \\
crystm02 & 1.40E+04 & 3.23E+05 & 1.66E-03 & 2.31E+01 & 1.76E-12 & 7.03E-15 & 4.49E+02 & 0.00E+00 & 3.98E-01 \\
finan512 & 7.48E+04 & 5.97E+05 & 1.07E-04 & 7.99E+00 & 2.82E+01 & 9.47E-01 & 8.50E+01 & 9.44E-01 & 6.43E-01 \\
bcsstk11 & 1.47E+03 & 3.42E+04 & 1.58E-02 & 2.32E+01 & 6.56E+08 & 2.96E+00 & 5.25E+08 & 0.00E+00 & 6.75E-04 \\
bloweybq & 1.00E+04 & 5.00E+04 & 5.00E-04 & 5.00E+00 & 5.00E+03 & 1.56E-15 & 4.17E+18 & 0.00E+00 & 1.00E-04 \\
nasa1824 & 1.82E+03 & 3.92E+04 & 1.18E-02 & 2.15E+01 & 2.12E+07 & -3.60E+00 & 1.42E+07 & 0.00E+00 & 7.93E-15 \\
mhd3200b & 3.20E+03 & 1.83E+04 & 1.79E-03 & 5.72E+00 & 2.20E+00 & 1.37E-13 & 2.02E+13 & 0.00E+00 & 1.98E-04 \\
msc10848 & 1.08E+04 & 1.23E+06 & 1.05E-02 & 1.13E+02 & 6.30E+11 & 6.32E+01 & 3.29E+10 & 1.84E-03 & 5.32E-03 \\
\bottomrule
\end{longtable}
\end{landscape}

\FloatBarrier \ \\ 
 \section*{Appendix - Mixed Models}\label{appendix:mixed_models}
\begin{table}
\centering
\caption{Probability of accuracy = 1 and slowdowns < 1.5}
\label{tab:MLM_merged}
\begin{tabular}{lrr}
\toprule
classifier & P(accuracy = 1) & P(slowdowns < 1.5) \\
\midrule
Random Forest\_256 & 0.60 & 0.71 \\
Random Forest\_32 & 0.56 & 0.68 \\
Random Forest\_64 & 0.56 & 0.68 \\
Random Forest\_128 & 0.55 & 0.68 \\
SVC\_256 & 0.53 & 0.61 \\
Gradient Boosting\_256 & 0.52 & 0.71 \\
SVC\_32 & 0.52 & 0.61 \\
K-Nearest Neighbors\_32 & 0.51 & 0.65 \\
SVC\_128 & 0.51 & 0.61 \\
SVC\_64 & 0.50 & 0.59 \\
Gradient Boosting\_32 & 0.48 & 0.68 \\
K-Nearest Neighbors\_64 & 0.47 & 0.61 \\
Gradient Boosting\_64 & 0.46 & 0.69 \\
Gradient Boosting\_128 & 0.46 & 0.70 \\
K-Nearest Neighbors\_256 & 0.44 & 0.65 \\
K-Nearest Neighbors\_128 & 0.42 & 0.61 \\
MLP Classifier\_32 & 0.41 & 0.71 \\
Decision Tree\_256 & 0.41 & 0.72 \\
Logistic Regression\_32 & 0.37 & 0.68 \\
Decision Tree\_32 & 0.37 & 0.70 \\
Decision Tree\_64 & 0.36 & 0.71 \\
Decision Tree\_128 & 0.36 & 0.70 \\
MLP Classifier\_64 & 0.34 & 0.68 \\
MLP Classifier\_128 & 0.32 & 0.72 \\
MLP Classifier\_256 & 0.31 & 0.74 \\
Logistic Regression\_64 & 0.29 & 0.72 \\
Gaussian Naive Bayes\_128 & 0.29 & 0.74 \\
Logistic Regression\_128 & 0.28 & 0.74 \\
Logistic Regression\_256 & 0.28 & 0.77 \\
Gaussian Naive Bayes\_32 & 0.22 & 0.68 \\
Gaussian Naive Bayes\_64 & 0.20 & 0.67 \\
Gaussian Naive Bayes\_256 & 0.20 & 0.72 \\
\bottomrule
\end{tabular}
\end{table}
\FloatBarrier \ \\


\begin{thebibliography}{10}
\expandafter\ifx\csname url\endcsname\relax
  \def\url#1{\texttt{#1}}\fi
\expandafter\ifx\csname urlprefix\endcsname\relax\def\urlprefix{URL }\fi
\expandafter\ifx\csname href\endcsname\relax
  \def\href#1#2{#2} \def\path#1{#1}\fi

\bibitem{balaprakash2018autotuning}
P.~Balaprakash, J.~Dongarra, T.~Gamblin, M.~Hall, J.~K. Hollingsworth, B.~Norris, R.~Vuduc, Autotuning in high-performance computing applications, Proceedings of the IEEE 106~(11) (2018) 2068--2083.
\newblock \href {https://doi.org/10.1109/JPROC.2018.2841200} {\path{doi:10.1109/JPROC.2018.2841200}}.

\bibitem{bauer2015quiet}
P.~Bauer, A.~Thorpe, G.~Brunet, \href{10.1038/nature14956}{The quiet revolution of numerical weather prediction}, Nature 525~(7567) (2015) 47--55.
\newline\urlprefix\url{10.1038/nature14956}

\bibitem{gasparini21HybridParallelIterativeSparseLinearSolverFramework}
L.~Gasparini, J.~R. Rodrigues, D.~A. Augusto, L.~M. Carvalho, C.~Conopoima, P.~Goldfeld, J.~Panetta, J.~P. Ramirez, M.~Souza, M.~O. Figueiredo, V.~M. Leite, Hybrid parallel iterative sparse linear solver framework for reservoir geomechanical and flow simulation, Journal of Computational Science 51 (2021) 101330.
\newblock \href {https://doi.org/10.1016/j.jocs.2021.101330} {\path{doi:10.1016/j.jocs.2021.101330}}.

\bibitem{gaganis2012machine}
\href{https://doi.org/10.2118/154505-MS}{Machine learning methods to speed up compositional reservoir simulation}, Vol. All Days of SPE Europec featured at EAGE Conference and Exhibition.
\newblock \href {http://arxiv.org/abs/https://onepetro.org/SPEEURO/proceedings-pdf/12EURO/All-12EURO/SPE-154505-MS/1612365/spe-154505-ms.pdf} {\path{arXiv:https://onepetro.org/SPEEURO/proceedings-pdf/12EURO/All-12EURO/SPE-154505-MS/1612365/spe-154505-ms.pdf}}, \href {https://doi.org/10.2118/154505-MS} {\path{doi:10.2118/154505-MS}}.
\newline\urlprefix\url{https://doi.org/10.2118/154505-MS}

\bibitem{scott2023introduction}
J.~Scott, M.~T{\r{u}}ma, Algorithms for sparse linear systems, Springer International Publishing, Cham, 2023.
\newblock \href {https://doi.org/10.1007/978-3-031-25820-6\_1} {\path{doi:10.1007/978-3-031-25820-6\_1}}.

\bibitem{saad2003iterative}
Y.~Saad, Iterative methods for sparse linear systems, SIAM, 2003.

\bibitem{benzi2002preconditioning}
M.~Benzi, Preconditioning techniques for large linear systems: a survey, Journal of Computational Physics 182~(2) (2002) 418--477.
\newblock \href {https://doi.org/10.1006/jcph.2002.7176} {\path{doi:10.1006/jcph.2002.7176}}.

\bibitem{bell2008efficient}
N.~Bell, M.~Garland, Efficient sparse matrix-vector multiplication on {CUDA}, Tech. rep., Nvidia Technical Report NVR-2008-004, Nvidia Corporation (2008).

\bibitem{Meijerink1977AnIS}
J.~A. Meijerink, H.~A. van~der Vorst, An iterative solution method for linear systems of which the coefficient matrix is a symmetric {M-matrix}, Mathematics of Computation 31 (1977) 148--162.

\bibitem{stuben2001AReview}
K.~Stüben, \href{https://www.sciencedirect.com/science/article\\/pii/B978044450617750015X}{A review of algebraic multigrid}, in: C.~Brezinski, L.~Wuytack (Eds.), Numerical Analysis: Historical Developments in the 20th Century, Elsevier, Amsterdam, 2001, pp. 331--359.
\newblock \href {https://doi.org/10.1016/B978-0-444-50617-7.50015-X} {\path{doi:10.1016/B978-0-444-50617-7.50015-X}}.
\newline\urlprefix\url{https://www.sciencedirect.com/science/article\\/pii/B978044450617750015X}

\bibitem{Davis2016ASO}
T.~A. Davis, S.~Rajamanickam, W.~M. Sid-Lakhdar, A survey of direct methods for sparse linear systems, Acta Numerica (2016) 383–566\href {https://doi.org/10.1017/S0962492916000076} {\path{doi:10.1017/S0962492916000076}}.

\bibitem{Li2015PerformanceAA}
K.~Li, W.~Yang, K.~Li, Performance analysis and optimization for {SpMV} on {GPU} using probabilistic modeling, IEEE Transactions on Parallel and Distributed Systems 26~(1) (2015) 196--205.
\newblock \href {https://doi.org/10.1109/TPDS.2014.2308221} {\path{doi:10.1109/TPDS.2014.2308221}}.

\bibitem{Zhang2011AQP}
Y.~Zhang, J.~D. Owens, A quantitative performance analysis model for {GPU} architectures, in: 2011 IEEE 17th International Symposium on High Performance Computer Architecture, 2011, pp. 382--393.
\newblock \href {https://doi.org/10.1109/HPCA.2011.5749745} {\path{doi:10.1109/HPCA.2011.5749745}}.

\bibitem{nisa2018effective}
I.~Nisa, C.~Siegel, A.~S. Rajam, A.~Vishnu, P.~Sadayappan, Effective machine learning based format selection and performance modeling for {SpMV} on {GPUs}, in: 2018 IEEE International Parallel and Distributed Processing Symposium Workshops (IPDPSW), IEEE, 2018, pp. 1056--1065.
\newblock \href {https://doi.org/10.1109/IPDPSW.2018.00164} {\path{doi:10.1109/IPDPSW.2018.00164}}.

\bibitem{Li2021AdaptiveSO}
M.~Li, Y.~Ao, C.~Yang, Adaptive {SpMV/SpMSpV} on {GPUs} for input vectors of varied sparsity, IEEE Transactions on Parallel and Distributed Systems 32~(7) (2021) 1842--1853.
\newblock \href {https://doi.org/10.1109/TPDS.2020.3040150} {\path{doi:10.1109/TPDS.2020.3040150}}.

\bibitem{barreda2020performance}
M.~Barreda, M.~F. Dolz, M.~A. Casta{\~n}o, P.~Alonso-Jord{\'a}, E.~S. Quintana-Orti, \href{doi.org/10.1007/s11227-020-03186-1}{Performance modeling of the sparse matrix--vector product via convolutional neural networks}, The Journal of Supercomputing 76~(11) (2020) 8883--8900.
\newline\urlprefix\url{doi.org/10.1007/s11227-020-03186-1}

\bibitem{yamada2018preconditioner}
K.~Yamada, T.~Katagiri, H.~Takizawa, K.~Minami, M.~Yokokawa, T.~Nagai, M.~Ogino, Preconditioner auto-tuning using deep learning for sparse iterative algorithms, in: 2018 Sixth International Symposium on Computing and Networking Workshops (CANDARW), IEEE, 2018, pp. 257--262.
\newblock \href {https://doi.org/10.1109/CANDARW.2018.00055} {\path{doi:10.1109/CANDARW.2018.00055}}.

\bibitem{li2022survey}
Z.~Li, F.~Liu, W.~Yang, S.~Peng, J.~Zhou, A survey of convolutional neural networks: analysis, applications, and prospects, IEEE Transactions on neural networks and learning systems 33~(12) (2022) 6999--7019.
\newblock \href {https://doi.org/10.1109/TNNLS.2021.3084827} {\path{doi:10.1109/TNNLS.2021.3084827}}.

\bibitem{falch2017machine}
T.~L. Falch, A.~C. Elster, \href{https://onlinelibrary.wiley.com/doi/abs/10.1002/cpe.4029}{Machine learning-based auto-tuning for enhanced performance portability of {OpenCL} applications}, Concurrency and Computation: Practice and Experience 29~(8) (2017) e4029, e4029 cpe.4029.
\newblock \href {https://doi.org/10.1002/cpe.4029} {\path{doi:10.1002/cpe.4029}}.
\newline\urlprefix\url{https://onlinelibrary.wiley.com/doi/abs/10.1002/cpe.4029}

\bibitem{tuncer2017diagnosing}
O.~Tuncer, E.~Ates, Y.~Zhang, A.~Turk, J.~Brandt, V.~J. Leung, M.~Egele, A.~K. Coskun, Diagnosing performance variations in {HPC} applications using machine learning, in: J.~M. Kunkel, R.~Yokota, P.~Balaji, D.~Keyes (Eds.), High Performance Computing: 32nd International Conference, ISC High Performance 2017, Frankfurt, Germany, June 18--22, 2017, Proceedings 32, Springer International Publishing, 2017, pp. 355--373.

\bibitem{memeti2019using}
S.~Memeti, S.~Pllana, A.~Binotto, J.~Ko{\l}odziej, I.~Brandic, \href{10.1007/s00607-018-0614-9}{Using meta-heuristics and machine learning for software optimization of parallel computing systems: a systematic literature review}, Computing 101~(8) (2019) 893--936.
\newline\urlprefix\url{10.1007/s00607-018-0614-9}

\bibitem{sedaghati2015automatic}
N.~Sedaghati, T.~Mu, L.-N. Pouchet, S.~Parthasarathy, P.~Sadayappan, Automatic selection of sparse matrix representation on {GPUs}, in: Proceedings of the 29th ACM on International Conference on Supercomputing, ICS '15, Association for Computing Machinery, 2015, pp. 99--108.
\newblock \href {https://doi.org/10.1145/2751205.2751244} {\path{doi:10.1145/2751205.2751244}}.

\bibitem{cui2016code}
H.~Cui, S.~Hirasawa, H.~Takizawa, H.~Kobayashi, A code selection mechanism using deep learning, in: 2016 IEEE 10th International Symposium on Embedded Multicore/Many-core Systems-on-Chip (MCSOC), IEEE, 2016, pp. 385--392.
\newblock \href {https://doi.org/10.1109/MCSoC.2016.46} {\path{doi:10.1109/MCSoC.2016.46}}.

\bibitem{peairs2011using}
L.~Peairs, T.-Y. Chen, Using reinforcement learning to vary the m in {GMRES(m)}, Procedia Computer Science 4 (2011) 2257--2266.
\newblock \href {https://doi.org/10.1016/j.procs.2011.04.246} {\path{doi:10.1016/j.procs.2011.04.246}}.

\bibitem{bhowmick2006application}
S.~Bhowmick, V.~Eijkhout, Y.~Freund, E.~Fuentes, D.~Keyes, Application of machine learning in selecting sparse linear solvers (2006).

\bibitem{dufrechou2019automatic}
E.~Dufrechou, P.~Ezzatti, E.~S. Quintana-Ort{\'\i}, Automatic selection of sparse triangular linear system solvers on {GPUs} through machine learning techniques, in: 2019 31st International Symposium on Computer Architecture and High Performance Computing (SBAC-PAD), IEEE, 2019, pp. 41--47.
\newblock \href {https://doi.org/10.1109/SBAC-PAD.2019.00020} {\path{doi:10.1109/SBAC-PAD.2019.00020}}.

\bibitem{funk2022prediction}
Y.~Funk, M.~G{\"o}tz, H.~Anzt, Prediction of optimal solvers for sparse linear systems using deep learning, in: Proceedings of the 2022 SIAM Conference on Parallel Processing for Scientific Computing, SIAM, 2022, pp. 14--24.
\newblock \href {https://doi.org/10.1137/1.9781611977141.2} {\path{doi:10.1137/1.9781611977141.2}}.

\bibitem{gotz2018machine}
M.~G{\"o}tz, H.~Anzt, Machine learning-aided numerical linear algebra: convolutional neural networks for the efficient preconditioner generation, in: 2018 IEEE/ACM 9th Workshop on Latest Advances in Scalable Algorithms for Large-Scale Systems (scalA), IEEE, 2018, pp. 49--56.
\newblock \href {https://doi.org/10.1109/ScalA.2018.00010} {\path{doi:10.1109/ScalA.2018.00010}}.

\bibitem{ackmann2020machine}
J.~Ackmann, P.~D. Düben, T.~N. Palmer, P.~K. Smolarkiewicz, Machine-learned preconditioners for linear solvers in geophysical fluid flows (2020).
\newblock \href {http://arxiv.org/abs/2010.02866} {\path{arXiv:2010.02866}}.

\bibitem{taghibakhshi2021optimization}
A.~Taghibakhshi, S.~MacLachlan, L.~Olson, M.~West, \href{https://proceedings.neurips.cc/paper\_files/paper/2021/file\\ /6531b32f8d02fece98ff36a64a7c8260-Paper.pdf}{Optimization-based algebraic multigrid coarsening using reinforcement learning}, in: M.~Ranzato, A.~Beygelzimer, Y.~Dauphin, P.~Liang, J.~W. Vaughan (Eds.), Advances in Neural Information Processing Systems, Vol.~34, Curran Associates, Inc., 2021, pp. 12129--12140.
\newline\urlprefix\url{https://proceedings.neurips.cc/paper\_files/paper/2021/file\\ /6531b32f8d02fece98ff36a64a7c8260-Paper.pdf}

\bibitem{saad1994ilut}
Y.~Saad, {ILUT}: A dual threshold incomplete {LU} factorization, Numerical linear algebra with applications 1~(4) (1994) 387--402.
\newblock \href {https://doi.org/10.1002/nla.1680010405} {\path{doi:10.1002/nla.1680010405}}.

\bibitem{kolodziej2019suitesparse}
S.~P. Kolodziej, M.~Aznaveh, M.~Bullock, J.~David, T.~A. Davis, M.~Henderson, Y.~Hu, R.~Sandstrom, The {SuiteSparse Matrix Collection} website interface, Journal of Open Source Software 4~(35) (2019) 1244.
\newblock \href {https://doi.org/10.21105/joss.01244} {\path{doi:10.21105/joss.01244}}.

\bibitem{ssgetpy}
S.~Raghunathan, \href{https://github.com/drdarshan/ssgetpy}{{SSGETPY}: Search and download sparse matrices from the {SuiteSparse Matrix Collection}}, {GitHub} repository (2023) [cited 2023-10-20T10:45:13].
\newline\urlprefix\url{https://github.com/drdarshan/ssgetpy}

\bibitem{matlab}
{MATLAB}, MathWorks, version 9.14.0.2254940 (R2023a) (2023).

\bibitem{ImagePrecGitHub}
M.~Souza, \href{https://github.com/michaelsouza/imageprec/blob/\\/main/README.md}{Dataset: Image and scalar-based approaches in preconditioner selection}, {GitHub} repository (2023) [cited 2023-10-26T18:20:10].
\newline\urlprefix\url{https://github.com/michaelsouza/imageprec/blob/\\/main/README.md}

\bibitem{guyon2006introduction}
I.~Guyon, A.~Elisseeff, An introduction to feature extraction, in: Feature extraction: foundations and applications, Springer Berlin Heidelberg, Berlin, Heidelberg, 2006, pp. 1--25.
\newblock \href {https://doi.org/10.1007/978-3-540-35488-8\_1} {\path{doi:10.1007/978-3-540-35488-8\_1}}.

\bibitem{petsc-web-page}
S.~Balay, S.~Abhyankar, M.~F. Adams, S.~Benson, J.~Brown, P.~Brune, K.~Buschelman, E.~Constantinescu, L.~Dalcin, A.~Dener, V.~Eijkhout, J.~Faibussowitsch, W.~D. Gropp, V.~Hapla, T.~Isaac, P.~Jolivet, D.~Karpeev, D.~Kaushik, M.~G. Knepley, F.~Kong, S.~Kruger, D.~A. May, L.~C. McInnes, R.~T. Mills, L.~Mitchell, T.~Munson, J.~E. Roman, K.~Rupp, P.~Sanan, J.~Sarich, B.~F. Smith, S.~Zampini, H.~Zhang, H.~Zhang, J.~Zhang, {PETSc/TAO} users manual, Tech. Rep. ANL-21/39 - Revision 3.20, Argonne National Laboratory (2023).
\newblock \href {https://doi.org/10.2172/1968587} {\path{doi:10.2172/1968587}}.

\bibitem{scikit-learn}
F.~Pedregosa, G.~Varoquaux, A.~Gramfort, V.~Michel, B.~Thirion, O.~Grisel, M.~Blondel, P.~Prettenhofer, R.~Weiss, V.~Dubourg, J.~Vanderplas, A.~Passos, D.~Cournapeau, M.~Brucher, M.~Perrot, {{\'E}}douard Duchesnay, \href{http://jmlr.org/papers/v12/pedregosa11a.html}{Scikit-learn: Machine learning in {Python}}, Journal of Machine Learning Research 12~(85) (2011) 2825--2830.
\newline\urlprefix\url{http://jmlr.org/papers/v12/pedregosa11a.html}

\bibitem{chollet2015keras}
F.~Chollet, et~al., \href{https://keras.io}{Keras} (2015) [cited 2023-09-01T11:27:26].
\newline\urlprefix\url{https://keras.io}

\bibitem{tensorflow2015-whitepaper}
M.~Abadi, A.~Agarwal, P.~Barham, E.~Brevdo, Z.~Chen, C.~Citro, G.~S. Corrado, A.~Davis, J.~Dean, M.~Devin, S.~Ghemawat, I.~Goodfellow, A.~Harp, G.~Irving, M.~Isard, Y.~Jia, R.~Jozefowicz, L.~Kaiser, M.~Kudlur, J.~Levenberg, D.~Man\'{e}, R.~Monga, S.~Moore, D.~Murray, C.~Olah, M.~Schuster, J.~Shlens, B.~Steiner, I.~Sutskever, K.~Talwar, P.~Tucker, V.~Vanhoucke, V.~Vasudevan, F.~Vi\'{e}gas, O.~Vinyals, P.~Warden, M.~Wattenberg, M.~Wicke, Y.~Yu, X.~Zheng, \href{https://www.tensorflow.org/}{{TensorFlow}: Large-scale machine learning on heterogeneous systems}, software available from tensorflow.org (2015) [cited 2023-09-01T10:34:13].
\newline\urlprefix\url{https://www.tensorflow.org/}

\bibitem{kotsiantis2007supervised}
S.~B. Kotsiantis, I.~Zaharakis, P.~Pintelas, et~al., Supervised machine learning: A review 

\bibitem{Rawat2017DeepCN}
W.~Rawat, Z.~Wang, Deep convolutional neural networks for image classification: a comprehensive review, Neural Computation 29~(9) (2017) 2352--2449.
\newblock \href {https://doi.org/10.1162/neco\_a\_00990} {\path{doi:10.1162/neco\_a\_00990}}.

\end{thebibliography}
\end{document}